\documentclass[11pt]{article}
\usepackage[utf8]{inputenc}
\usepackage{fullwidth}
\usepackage{color}
\usepackage[top=1.5in, bottom=1.5in, left=1.2in, right=1.2in]{geometry}
\usepackage{url}
\usepackage{amsmath}
\usepackage{amstext}
\usepackage{graphicx}
\usepackage[raggedright]{titlesec}
\usepackage[font=small,labelfont=bf]{caption}
\usepackage{subfig}

\title{Algebraic Structure and Complexity of Bootstrap Percolation with External Inputs}

\author{Saptarshi Pal, University of Waterloo \and Chrystopher L. Nehaniv, Waterloo Algebraic Intelligence and \\ Computation Laboratory, University of Waterloo}

\begin{document}
\maketitle 
\abstract{ In this paper a modification of the standard Bootstrap Percolation model is introduced. In our modification a discrete time update rule is constructed that allows for non-monotonicity - unlike its classical counterpart. External inputs to drive the system into desirable states are also included in the model. The algebraic structure and complexity properties of the system are inferred by studying the system's holonomy decomposition. We introduce methods of inferring the pools of reversibility for the system. Dependence of system complexity on process parameters is presented and discussed.
}

\section{Introduction}
\label{sec:1}
Bootstrap percolation is a process studied in statistical mechanics where cells in a lattice or any other space (for example, nodes in a random graph) exist in binary states --- active or inactive. Given any initial configuration of states, the states of the nodes evolve with discrete time based on some predefined update rule. Bootstrap percolation is an example of a boolean network. The most popular update rule of bootstrap percolation is parameterized with a threshold $k$ as follows:

 \[
    x^i(t+1) = \begin{cases}
        x^i(t) & \quad \text{if } x^i(t) = 1\\
        1 & \quad \text{if $x^i(t)$ = 0 and  } \sum_{j \in N(i)}x^i(t) \geq k\\
        0 & \quad \text{if $x^i(t)$ = 0 and  } \sum_{j \in N(i)}x^i(t) < k,
        \end{cases}
  \]
  
where $N(i)$ is the set of neighbours of node $i$ and $x^i(t)$ is the boolean state of the node $i$ at discrete time $t$: $x^i(t) = 1$ representing active state and $x^i(t) = 0$ representing the inactive state for node $i$ at time $t$. The standard interpretation of the states may be inverted in some cases causing no loss of generality. There can be other update rules of a bootstrap percolation process like in \cite{einarsson2014bootstrap} but all variants have the same property of being homogeneous and local. The standard model of bootstrap percolation is monotone with evolution of discrete time. Some variants of the bootstrap percolation, with the introduction of excitatory nodes and inhibitory nodes \cite{einarsson2014bootstrap}, display non-monotone behavior. The model of bootstrap percolation with inhibition has been useful in studying the phenomenon of input normalization in neurons \cite{carandini2013erratum}. In the next section a modification to the standard bootstrap percolation is introduced that also exhibits a non-monotonic behavior. 
Bootstrap percolation models have also been used to study impacts of external perturbations to models of weighted trade networks \cite{fan2014state}. 
In \cite{fan2014state}, each node - a nation state - is modelled to be in a binary state. A node in an active state represents that the corresponding nation state is in a normal state, i.e. it has imports and exports in a stable level. An abnormal or inactive state for a nation indicates that it is facing a trading disaster and it is leading towards a volatile economy. Since an abnormal country can influence the state of its neighbours by potentially turning them into abnormals, bootstrap percolation is used for modeling the cascading reaction of any initial perturbation. 
Despite being an interesting model of the spread of trading and economic disasters in a trade network, certain major assumptions of the model make it less flexible to real settings. The model does not allow nodes going from inactive to active states through a local process in the same way they go from active to inactive. Moreover the model, in general does not allow independent relapse or recovery of nodes - i.e without the assistance of local homogenous update rule. Even though our modified bootstrap percolation model was motivated by drawbacks and inadequacies in the fault propagation model in trade-networks, one can use the structural framework of our model to represent other non-monotone local-homogenous processes. 

\section{The modified Bootstrap Percolation Model with external inputs}

The new model description can be broken down into two independent parts:
\begin{enumerate}
    \item The non-monotone bootstrap percolation process.
    \item External inputs forcing certain network states to certain target network states.
\end{enumerate}
The second part of the model has a superiority over the first part. After the formal description of these parts it will be clear as to how their hierarchy is defined.\\ \\
\textbf{The non-Monotone Bootstrap Percolation Process:}\\ \\
This part of the model corresponds to the discrete step update rule. As opposed to the standard bootstrap percolation model, this modification has two process parameters $k_1$, $k_2$ instead of a single one $k$. The process is defined as: \\
\[ x^i(t+1) = \begin{cases} 
     1, & \quad \text{if $x^i(t)$ = 0 and  } \sum_{j \in N(i)}x^i(t) \geq k_1\\
        0, & \quad \text{if $x^i(t)$ = 0 and  } \sum_{j \in N(i)}x^i(t) < k_1\\ \\
        0, & \quad \text{if $x^i(t)$ = 1 and  } \sum_{j \in N(i)}x^i(t) \leq |N(i)| - k_2 \\
        1, & \quad \text{if $x^i(t)$ = 1 and  } \sum_{j \in N(i)}x^i(t) > |N(i)| - k_2
   \end{cases}
\]

In brief, the above process can be described textually as: (a) If a node is inactive and $k_1$ or more of its neighbours are active, it will turn active. (b) If a node is active and $k_2$ or more neighbours are inactive it will become inactive.  Otherwise the state of nodes remains unchanged.\\
Note that now, with these changes we have an identical jump back rule from active to inactive like jumping from inactive to active. Unless there are any external inputs (the second part of the model), these update rules are obeyed at every discrete time step to evaluate the state of the network at the next time step. This part forms the local-homogenous update rule for our model.\\ \\
\noindent
\textbf{External Inputs:} \\ \\
Only the above modification does not allow nodes to have independent transitions to active or inactive states, i.e without the assistance of the homogenous update rule. In order to account for that, external forced inputs are introduced. Before defining the external forced inputs let's define what is meant by state of the network. Since a node can exist is two states and there are $N$ nodes in the network, there can be $2^N$ possibilities of network state. Network state at discrete time step $t$ is defined as $X(t) = (x^1(t),x^2(t),...,x^N(t))$. The network state is an element from $\{0,1\}^N$. For convenience, each network state will be encoded with its decimal equivalent added with 1. So, for example the state $(1,0,0,...,0)$ represents that only the node labelled $1$ is active ($1$) and rest are inactive ($0$). If there were 4 nodes in the network the state $(1,0,0,0)$ would have been encoded as 9. The $2^N$ states of a network are represented as $X^1$, $X^2$, ..., $X^{2^N}$\\
The external inputs are analogous to constant maps except for a crucial distinction: that some states of the network are non-forceable. This brings us to the third parameter of our model which is the set of states of the network which cannot be forced into a target set by any constant map. Formally, this part of the model can be expressed as:
\\For a network of size $N$ the set of all possible external inputs is:
$$S_{\text{all}} = \{s_1, s_2, ..., s_{2^N} \}$$
Let's define a one to one function $F$ as:
$F: s_i \mapsto X^i$.
If $S$ is defined as the set of non-forceable  states of the network, then if at time step $t$ there is an external input $s_k$ to the process,
$$X(t+1) = F(s_k) \quad \text{iff } X(t) \notin S$$
If $X(t) \in S$ then $X(t + 1)$ is determined by the local update rule of the process. At any time step $t$ the model does not permit more than one external input to the process. \\
This part of the model has a hierarchy over the previous part. If there is any external input to the process, the external input part takes precedence over the local homogenous update rule.

\section{Holonomy Decomposition and Modeling the Process as a Transformation Semigroup}
\begin{figure}
    \centering
    \includegraphics[width = 0.5\linewidth]{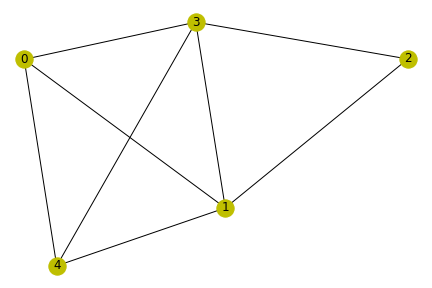}
    \caption{Random Graph with 5 nodes chosen for this study. A total of 32 states starting from 1 till 32 is possible for this network, for all the different boolean labellings of the nodes. Numbers inside the node indicate their index in the network state.}
    \label{Fig 1}
\end{figure}

In this paper the complexity and behaviour of such a modified model is analyzed by observing its algebraic structure. There are some obvious questions this paper aims to investigate. Since the process is non-monotone it is of interest to investigate whether a particular instance (a set of parameters $k_1$, $k_2$ and $S$) of this model falls into cycles or dies down into a single state. It may also be of interest to know whether it is possible to trigger the system into cycle - if yes, what is the fastest possible way of achieving that. Intuitively, the complexity of the process increases if the cardinality of $S$ increases. In  this section two measures are defined to quantify the algebraic complexity of the process. 
We will be using the method of computational Holonomy Decomposition of Transformation Actions to analyze this model \cite{egri2015computational}. The computational holonomy decomposition of any discrete-time, discrete state process requires the process to be represented as a transformation action. Before defining a transformation action representation the modified Bootstrap percolation model let's look at holonomy decomposition and the advantages of using it to analyze our system.\\
A transformation action (also called a transformation semigroup) is defined as a set of functions that maps a set to itself and is closed under function composition. It is represented as $(A,S)$ where $S$ is the set of functions acting on the set $A$. (This $S$ is not to be confused with the set of non-forceable states as defined in the model above.) Let $a\cdot s$ denote resulting state $s(a)$ in $A$ resulting from applying transformation $s\in S$ to state $a\in A$. Similarly for $P\subseteq A$, $P\cdot s =\{a \cdot s \in A\,\vert \, a\in P\}$.   The \textit{extended image set} $I^*$ of the transformation action $(A,S)$ is defined as:
$$I^* = \{A\cdot s\,\vert\, s \in S\} \cup \{A\} \cup \{\{a\}\,\vert\, a \in A\},$$
A reflexive and transitive relation called the \textit{subduction relation} is defined on $I^*$ as:
$$P \leq_S Q \iff P=Q  \mbox{ or } \exists s \in S\text{ such that } P\subseteq Q\cdot s \text{ for } P,Q \in I^*$$
\textit{Height} of a singleton member $Q$ in $I^*$ is defined to be 0. For any other member $Q \in I^*$, the height is $i$ is defined as the length of the longest strict subduction chain to $Q$ in $I^*$ that ends in a singleton, the $Q$'s holonomy group $(B_Q, H_Q)$ can be defined as is done in $\cite{egri2005cycle}$. A permutation reset semigroup $(\Gamma_i, \overline{\Theta_i})$ for a height $i$, $i = 1,...,h$ is defined as the direct product of holonomy groups for subduction equivalence class representatives $Q$ in $I^*$ that have height $i$, augmented with the constant maps.\\
The \textit{holonomy decomposition theorem} \cite{egri2005cycle} states that any finite transformation semigroup \textit{divides} or is \textit{emulated by} a wreath product of its holonomy permutation reset transformation semigroups. Since holonomy decomposition is used to study the modified bootstrap percolation - the process is modeled as a transformation semigroup as follows:\\
\begin{equation}
(X = \{0,1\}^N, \langle (S_{\text{all}} \setminus \{F^{-1}(s) | s \in S\}) \cup \{t\}\rangle)
\end{equation}
Note that here $S$ is the set of non-forceable states and $S_{\text{all}}$ is the set of all external inputs. The angular brackets indicate that the semigroup acting on the state set $X$ is generated by association of elements of the set inside the brackets. The actions of the transformation action is defined as:
$$x \cdot t = y \implies X(t+1) = y \text{ if } X(t) = x, \quad x,y \in X \text{ and no external input on system.} $$
$$\text{ and } x\cdot s_i = F(s_i) \quad \text{ if } x \notin S$$
$$x \cdot s_i = x \cdot t \quad \text{ if } x \in S$$
For every parameter set $k_1,k_2,S$ of the model for a given graph, a transformation semigroup can be written for the process using Equation 1. For the remainder of this paper we show the results for the graph indicated in Fig \ref{Fig 1}.
\section{Methods, Results and Discussion}
All the analysis of holonomy decomposition, algebraic structure, hierarchy and complexity were carried out using the computational tools of GAP (Groups, Algorithm, Programming)- a system for computational discrete algebra \cite{GAP4} and one of its packages: SgpDec (package for semigroup decomposition) \cite{AttilaSgp}.
In this paper, using GAP/SgpDec we study the holonomy decomposition of finite transformation actions \cite{eilenberg1976,egri2015computational} to gains insights about the discrete system dynamics. Another concept that will be frequently used in this study is algebraic complexity of a process. We use two measures of complexity: (a) Krohn-Rhodes (KR) Complexity Measure (b) Aperiodic Complexity Measure. The KR complexity measure is a unique maximal hierarchical complexity measure satisfying the complexity axioms defined in \cite{nehaniv2000evolution}. It is defined as the the smallest number of permutation levels needed in any Krohn-Rhodes decomposition \cite{krohn1965algebraic}. In this paper we employ a computable upper bound on KR complexity of a transformation semigroup (TS) given as  the total number of levels of the holonomy decomposition of TS with groups in them. We also define the \textit{height complexity}  of a transformation semigroup  as the smallest number of levels needed in the holonomy decomposition of the transformation semigroup. It is to be noted that all the complexity measures reported in this paper are upper bounds of the actual KR complexity measure as GAP/SgpDec is never always guaranteed to produce a shortest Krohn-Rhodes decomposition for a transformation semigroup.\\
Using SgpDec it is possible to observe and study the skeleton of a transformation semigroup \cite{egri2010skeleton}. A skeleton of a transformation action is a pre-ordered structure on $I^*$ which encodes information about the different ways of traversing from the entire state set to singletons of the system. It is possible to evaluate the upper bound of the number of consecutive irreversible transformations required (from any initialization) to reach a death state (or terminal generalized limit cycle, i.e., permutation group) for the system from the skeleton of a transformation action.\\
Since there are a lot of parameter combinations possible for our model, we focus on three scenarios of the process model which are somewhat realistic in the context of trade networks to demonstrate how holonomy decomposition can be used to study the bootstrap percolation process. Using those parameter values we will observe how the use of holonomy decomposition helps to understand the model process in more depth. \\ 
\noindent
We study the following three scenarios of the modified Bootstrap percolation model. For all the three scenarios the values of $k_1$ and $k_2$ range from 1 to 5:
\begin{enumerate}
    \item \textbf{Scenario 1}: In this scenario the set $S$ (the set of non-forceable states) is  empty.
    \item \textbf{Scenario 2}: In this scenario the set $S$ is taken as \{11,19,20,25\}.
    \item \textbf{Scenario 3}: In this scenario the set $S$ is taken as \{8,15,20,22,29\}.
\end{enumerate}
The choices for the set $S$ in the above scenarios are not random. Their construction is dependent on the subject graph that has been chosen for this study. In the second scenario the motivation was to bunch five states into the non-forceable set $S$ which allow only two nodes to be active with the constraint that their cumulative degree does not exceed $6$. Similarly in the third scenario the degree constraint was upper bounded to $10$. The construction of these scenarios are realistic because in most applied cases, like the case of a trade network, making states non-forceable incurs certain costs.\\

\begin{figure}[h!]
\centering
\subfloat[]
{
    \includegraphics[width = 0.5\linewidth]{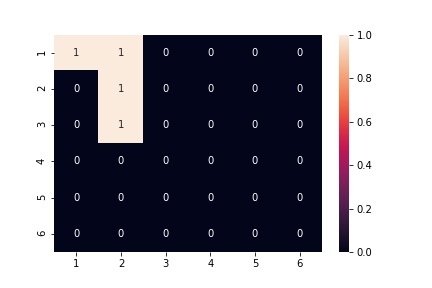}
    \label{scenario1a}
}
\subfloat[]
{
    \includegraphics[width = 0.5\linewidth]{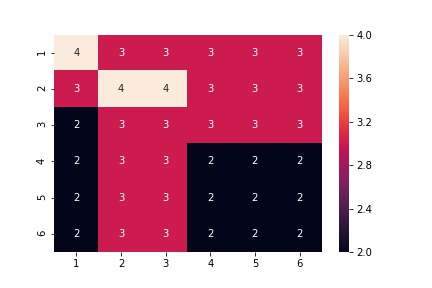}
    \label{scenario1b}
}
\caption{Complexity Analysis for Scenario 1. (a) Heat map for Krohn-Rhodes in scenario 1 (b) Heat map for  height complexity in scenario 1. The horizontal axis of the heat map corresponds to the parameter $k_1$. The vertical axis corresponds to the parameter $k_2$. The numbers inside the cells represent the complexity value at the corresponding values of $k_1$ and $k_2$.}
\end{figure}
\noindent
\textbf{Scenario 1 results:} The upper bound for complexity measures for scenario 1 are presented in Figure \ref{scenario1a} and \ref{scenario1b} with a heat map. For values of $(k_1,k_2) = (1,1),(2,1),(2,2)$ and $(2,3)$ the highest value for Krohn-Rhodes Complexity, i.e 1, is achieved.  For $k_1 > 3$ and $k_2 > 3$ the Krohn-Rhodes complexity of the system is zero. This implies that in that subspace of the parameter space it is possible to construct an embedding of the model system by using only banks of flip-flops as there exists no pools of reversibility in the system in them. The system never falls into cycles and always dies into a state. The holonomy components of the system in scenario 1 as found by SgpDec for $k_1, k_2 = (1,1)$ are:
Level 1: 19,\,
Level 2: 6,\,
Level 3: 2, \,
Level 4: (7,C4) \\
\begin{figure}[h!]
    \centering
    \includegraphics[width = 0.8\linewidth]{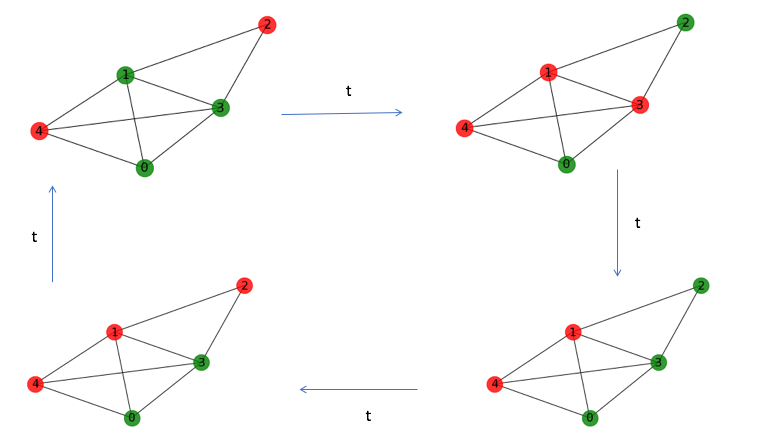}
    \caption{The C4 cycle in last level of the Holonomy Decomposition of Scenario 1 set at $(k_1,k_2) = (1,1)$. Red nodes indicate inactive nodes and green nodes indicate active nodes.}
    \label{Fig 4}
\end{figure}
\noindent
This encodes that the first level of the decomposition has 19 elements of extended image set being moved around by constant maps and identity. These transformation actions are called \textit{identity resets}. We also see identity resets in levels 2 and 3. In the fourth level of the decomposition we have a permutation group (7,C4) with a cyclic group $C4$ of order 4 acting on 7 elements, augmented with constant maps.\\
Although not presented in this report, it was found that the standard monotone bootstrap percolation model with constant maps (a simpler version than this one) has the simplest algebraic decompositions as it can be entirely built up with banks of flip flops. In other words the standard bootstrap percolation model with constant maps has Krohn-Rhodes complexity zero for all pairs of $(k_1, k_2)$ and any other modification of that model can be viewed as being an added complexity to the fundamental model. The case of $k_1, k_2 = (1,1)$ does however have cycles in it (indicated by permutation group in last level of decompostion). When the (7,C4) holonomy group is investigated further using SgpDec, it was found that the states \{6,12,10,14\} are being moved by the generator $t$ in a cycle, isomorphic to C4. In Fig~\ref{Fig 4} we see the physical interpretation of this. The cycle in which the graph falls into has been identified exactly. If the system falls into any one of the states shown in Fig \ref{Fig 4} it will forever be stuck there (since we are in scenario 1 there is no forced exit from the loop as well).
\begin{figure}[t]
\centering
\subfloat[]
{
    \includegraphics[width = 0.5\linewidth]{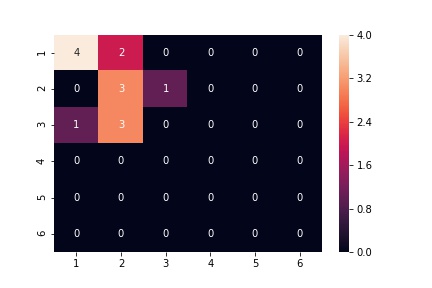}
    \label{fig5a}
}
\subfloat[]
{
    \includegraphics[width = 0.5\linewidth]{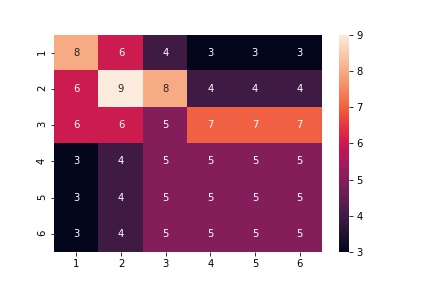}
    \label{fig5b}
} \\
\subfloat[]
{
    \includegraphics[width = 0.5\linewidth]{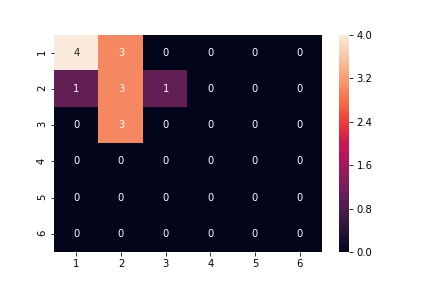}
    \label{fig6a}
}
\subfloat[]
{
    \includegraphics[width = 0.5\linewidth]{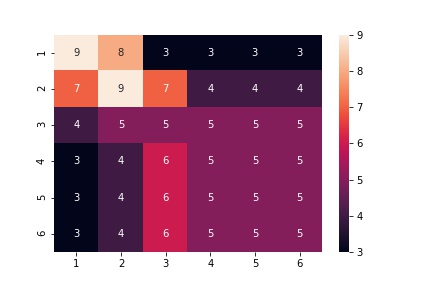}
    \label{fig6b}
}
\caption{Heat Map for complexity measures over parameters $k_1, k_2$ for scenario 2 and scenario 3 (a) KR Complexity Upper-bound Measure for scenario 2 (b) Height Complexity Measure for scenario 2. (c)  KR Complexity Upper-bound Measure for scenario 3 (d)  Height Complexity Measure for scenario 3. Highest KR complexity upper bound determined by SgpDec for scenarios 2 and 3 are both 4. }
\label{fig 5}
\end{figure}

\begin{figure}[t!]
    \centering
    \includegraphics[width = 0.8\linewidth]{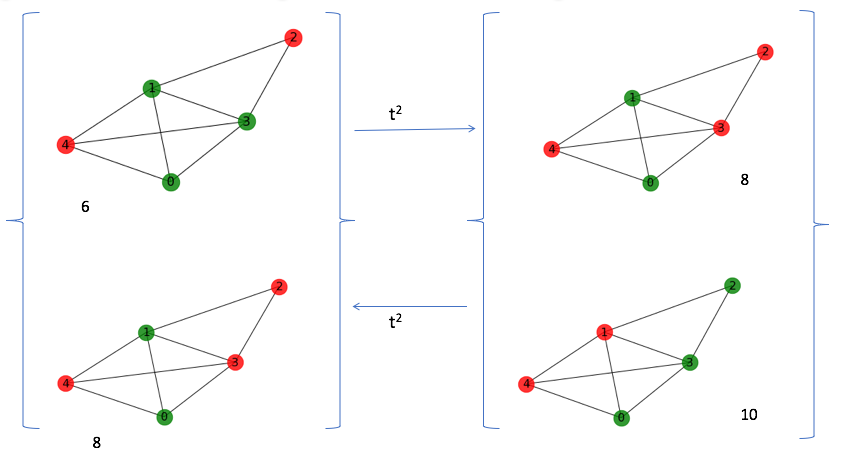}
    \caption{The C2 cycle in the $7^{th}$ level of the holonomy decomposition of Scenario 2 set at $(k_1,k_2) = (1,1)$. Red nodes are inactive nodes and green nodes are active nodes. The set \{6,8\} is mapped to \{8,10\} and back by $t^2$.}
    \label{c2}
\end{figure}

\noindent
Complexity increases in \textbf{Scenario 2} and \textbf{3} due to introduction of a non-forceable state as is observable in the heat maps of complexity in Fig \ref{fig 5}. Many more holonomy groups start to appear in the decomposition. Unlike scenario 1, in this scenario we have holonomy groups in levels in which there were none in the previous scenario. This implies that there are elements in the semigroup that move around sets of states together. To illustrate let's take an example. At $(k_1,k_2) = (1,1)$ for scenario 2, we find that the holonomy group (5,C2) appears in the penultimate level of decomposition (the $7^{th}$). On investigating this group closely we see that the set of states \{6,8\} was being moved around by $t^2$ to \{8,10\} and back. This is shown in Fig \ref{c2}. In this case since an element of the transformation semigroup (here $t^2$ meaning two clock ticks) moves set of states around in cycles, it is less intuitive to understand the actual cycles from the holonomy group itself. Looking at Natural Subsystems can help in having a better insight in these cases \cite{nehaniv2015symmetry}.\\
We also see that the highest KR complexity bound attained in Scenario 2 and 3 is 4. Scenario 3 is only a bit more complex than Scenario 2 in the height complexity measure. For some parameter tuples $(k_1,k_2)$ the upper bound of the number of irreversible actions required on the original state set to bring system to a death state is higher for Scenario 3 when compared to Scenario 2 (and of course, scenario 1). Note that in Scenario 2 and 3, just like Scenario 1, there exist threshold values for $k_1$ and $k_2$ above which KR complexity measure is zero as no pools of reversibility can exist in those parameter spaces. 

\section{Conclusion}

In this paper we introduced a modified non-monotone version of the classical bootstrap percolation that allowed external inputs and modelled it as an automaton in order to analyze its algebraic structure and corresponding complexity. Existence of cycles can be beneficial or detrimental depending on the system being studied. So their identification in a discrete event discrete time system process can be of interest. The complexity measures identified in this paper provide insight on that matter. A non-zero KR complexity measure for discrete systems indicates the existence of
subsets of the state space that can be moved around in cycles by sequence of combination of external operations. We also discuss analysis by holonomy decomposition of our bootstrap percolation model and discuss how it reveals hidden structures and information about the system. The decomposition reveals that our system's hierarchical construction can be emulated by building blocks of cyclic groups, simple non-abelian groups and flip-flops. Apart from cycles, other interesting system properties are also identified for our system by studying metrics like height complexity of a holonomy decomposition.
The non-monotone version of the Bootstrap percolation model under discussion shows higher complexity in structure than the classical Bootstrap Percolation model. Adding a non-forceable state set on top of that results in further increases in the complexity measures.
\section{Acknowledgements}
This work was supported in part by a Natural Sciences and Engineering Re-search Council of Canada (NSERC) grant, funding ref. RGPIN-2019-04669.
\bibliographystyle{abbrv}
\bibliography{references.bib}

\begin{thebibliography}{10}

\bibitem{carandini2013erratum}
M.~Carandini and D.~J. Heeger.
\newblock Erratum: Normalization as a canonical neural computation ({Nature
  Reviews Neuroscience} (2012) 13 (51-62)).
\newblock {\em Nature Reviews Neuroscience}, 14(2):152, 2013.

\bibitem{AttilaSgp}
A.~Egri-Nagy, C.~Nehaniv, and J.~D. Mitchell.
\newblock {SgpDec - software package for Hierarchical Composition and
  Decomposition of Permutation Groups and Transformation Semigroups}.
\newblock \url{https://github.com/gap-system/sgpdec}, 2015.

\bibitem{egri2005cycle}
A.~Egri-Nagy and C.~L. Nehaniv.
\newblock Cycle structure in automata and the holonomy decomposition.
\newblock {\em Acta Cybernetica}, 17(2):199--211, 2005.

\bibitem{egri2010skeleton}
A.~Egri-Nagy and C.~L. Nehaniv.
\newblock On the skeleton of a finite transformation semigroup.
\newblock In {\em Annales Mathematicae et Informaticae}, volume~37, pages
  77--84, 2010.

\bibitem{egri2015computational}
A.~Egri-Nagy and C.~L. Nehaniv.
\newblock Computational holonomy decomposition of transformation semigroups.
\newblock {\em arXiv preprint arXiv:1508.06345}, 2015.

\bibitem{eilenberg1976}
S.~Eilenberg.
\newblock {\em Automata, Languages and Machines}, volume~B.
\newblock Academic Press, 1976.

\bibitem{einarsson2014bootstrap}
H.~Einarsson, J.~Lengler, K.~Panagiotou, F.~Mousset, and A.~Steger.
\newblock Bootstrap percolation with inhibition.
\newblock {\em arXiv preprint arXiv:1410.3291}, 2014.

\bibitem{fan2014state}
Y.~Fan, S.~Ren, H.~Cai, and X.~Cui.
\newblock The state's role and position in international trade: A complex
  network perspective.
\newblock {\em Economic Modelling}, 39:71--81, 2014.

\bibitem{GAP4}
The GAP~Group.
\newblock {\em {GAP -- Groups, Algorithms, and Programming, Version 4.11.0}},
  2020.

\bibitem{krohn1965algebraic}
K.~Krohn and J.~Rhodes.
\newblock Algebraic theory of machines. i. prime decomposition theorem for
  finite semigroups and machines.
\newblock {\em Transactions of the American Mathematical Society},
  116:450--464, 1965.

\bibitem{nehaniv2015symmetry}
C.~L. Nehaniv, J.~Rhodes, A.~Egri-Nagy, P.~Dini, E.~R. Morris, G.~Horv{\'a}th,
  F.~Karimi, D.~Schreckling, and M.~J. Schilstra.
\newblock Symmetry structure in discrete models of biochemical systems: natural
  subsystems and the weak control hierarchy in a new model of computation
  driven by interactions.
\newblock {\em Philosophical Transactions of the Royal Society A: Mathematical,
  Physical and Engineering Sciences}, 373(2046):20140223, 2015.

\bibitem{nehaniv2000evolution}
C.~L. Nehaniv and J.~L. Rhodes.
\newblock The evolution and understanding of hierarchical complexity in biology
  from an algebraic perspective.
\newblock {\em Artificial Life}, 6(1):45--67, 2000.

\end{thebibliography}
\end{document}